\expandafter\chardef\csname pre amssym.def at\endcsname=\the\catcode`\@ 
\catcode`\@=11 
 
\def\undefine#1{\let#1\undefined} 
\def\newsymbol#1#2#3#4#5{\let\next@\relax 
 \ifnum#2=\@ne\let\next@\msafam@\else 
 \ifnum#2=\tw@\let\next@\msbfam@\fi\fi 
 \mathchardef#1="#3\next@#4#5} 
\def\mathhexbox@#1#2#3{\relax 
 \ifmmode\mathpalette{}{\m@th\mathchar"#1#2#3}%
 \else\leavevmode\hbox{$\m@th\mathchar"#1#2#3$}\fi} 
\def\hexnumber@#1{\ifcase#1 0\or 1\or 2\or 3\or 4\or 5\or 6\or 7\or 8\or 
 9\or A\or B\or C\or D\or E\or F\fi} 
 
\font\tenmsa=msam10 
\font\sevenmsa=msam7 
\font\fivemsa=msam5 
\newfam\msafam 
\textfont\msafam=\tenmsa 
\scriptfont\msafam=\sevenmsa 
\scriptscriptfont\msafam=\fivemsa 
\edef\msafam@{\hexnumber@\msafam} 
\mathchardef\dabar@"0\msafam@39 
\def\dashrightarrow{\mathrel{\dabar@\dabar@\mathchar"0\msafam@4B}} 
\def\dashleftarrow{\mathrel{\mathchar"0\msafam@4C\dabar@\dabar@}} 
 
\def\ulcorner{\delimiter"4\msafam@70\msafam@70 } 
\def\urcorner{\delimiter"5\msafam@71\msafam@71 } 
\def\llcorner{\delimiter"4\msafam@78\msafam@78 } 
\def\lrcorner{\delimiter"5\msafam@79\msafam@79 } 
\def\yen{{\mathhexbox@\msafam@55 }} 
\def\checkmark{{\mathhexbox@\msafam@58 }} 
\def\circledR{{\mathhexbox@\msafam@72 }} 
\def\maltese{{\mathhexbox@\msafam@7A }} 
 
\font\tenmsb=msbm10 
\font\sevenmsb=msbm7 
\font\fivemsb=msbm5 
\newfam\msbfam 
\textfont\msbfam=\tenmsb 
\scriptfont\msbfam=\sevenmsb 
\scriptscriptfont\msbfam=\fivemsb 
\edef\msbfam@{\hexnumber@\msbfam}

\catcode`\@=\csname pre amssym.def at\endcsname 
 
\expandafter\ifx\csname pre amssym.tex at\endcsname\relax \else \endinput\fi 
\expandafter\chardef\csname pre amssym.tex at\endcsname=\the\catcode`\@ 
\catcode`\@=11 
\newsymbol\boxdot 1200 
\newsymbol\boxplus 1201 
\newsymbol\boxtimes 1202 
\newsymbol\square 1003 
\newsymbol\blacksquare 1004 
\newsymbol\centerdot 1205 
\newsymbol\lozenge 1006 
\newsymbol\blacklozenge 1007 
\newsymbol\circlearrowright 1308 
\newsymbol\circlearrowleft 1309 
\undefine\rightleftharpoons 
\newsymbol\rightleftharpoons 130A 
\newsymbol\leftrightharpoons 130B 
\newsymbol\boxminus 120C 
\newsymbol\Vdash 130D 
\newsymbol\Vvdash 130E 
\newsymbol\vDash 130F 
\newsymbol\twoheadrightarrow 1310 
\newsymbol\twoheadleftarrow 1311 
\newsymbol\leftleftarrows 1312 
\newsymbol\rightrightarrows 1313 
\newsymbol\upuparrows 1314 
\newsymbol\downdownarrows 1315 
\newsymbol\upharpoonright 1316 
  
\newsymbol\downharpoonright 1317 
\newsymbol\upharpoonleft 1318 
\newsymbol\downharpoonleft 1319 
\newsymbol\rightarrowtail 131A 
\newsymbol\leftarrowtail 131B 
\newsymbol\leftrightarrows 131C 
\newsymbol\rightleftarrows 131D 
\newsymbol\Lsh 131E 
\newsymbol\Rsh 131F 
\newsymbol\rightsquigarrow 1320 
\newsymbol\leftrightsquigarrow 1321 
\newsymbol\looparrowleft 1322 
\newsymbol\looparrowright 1323 
\newsymbol\circeq 1324 
\newsymbol\succsim 1325 
\newsymbol\gtrsim 1326 
\newsymbol\gtrapprox 1327 
\newsymbol\multimap 1328 
\newsymbol\therefore 1329 
\newsymbol\because 132A 
\newsymbol\doteqdot 132B 
  
\newsymbol\triangleq 132C 
\newsymbol\precsim 132D 
\newsymbol\lesssim 132E 
\newsymbol\lessapprox 132F 
\newsymbol\eqslantless 1330 
\newsymbol\eqslantgtr 1331 
\newsymbol\curlyeqprec 1332 
\newsymbol\curlyeqsucc 1333 
\newsymbol\preccurlyeq 1334 
\newsymbol\leqq 1335 
\newsymbol\leqslant 1336 
\newsymbol\lessgtr 1337 
\newsymbol\backprime 1038 
\newsymbol\risingdotseq 133A 
\newsymbol\fallingdotseq 133B 
\newsymbol\succcurlyeq 133C 
\newsymbol\geqq 133D 
\newsymbol\geqslant 133E 
\newsymbol\gtrless 133F 
\newsymbol\sqsubset 1340 
\newsymbol\sqsupset 1341 
\newsymbol\vartriangleright 1342 
\newsymbol\vartriangleleft 1343 
\newsymbol\trianglerighteq 1344 
\newsymbol\trianglelefteq 1345 
\newsymbol\bigstar 1046 
\newsymbol\between 1347 
\newsymbol\blacktriangledown 1048 
\newsymbol\blacktriangleright 1349 
\newsymbol\blacktriangleleft 134A 
\newsymbol\vartriangle 134D 
\newsymbol\blacktriangle 104E 
\newsymbol\triangledown 104F 
\newsymbol\eqcirc 1350 
\newsymbol\lesseqgtr 1351 
\newsymbol\gtreqless 1352 
\newsymbol\lesseqqgtr 1353 
\newsymbol\gtreqqless 1354 
\newsymbol\Rrightarrow 1356 
\newsymbol\Lleftarrow 1357 
\newsymbol\veebar 1259 
\newsymbol\barwedge 125A 
\newsymbol\doublebarwedge 125B 
\undefine\angle 
\newsymbol\angle 105C 
\newsymbol\measuredangle 105D 
\newsymbol\sphericalangle 105E 
\newsymbol\varpropto 135F 
\newsymbol\smallsmile 1360 
\newsymbol\smallfrown 1361 
\newsymbol\Subset 1362 
\newsymbol\Supset 1363 
\newsymbol\Cup 1264 
  
\newsymbol\Cap 1265 
  
\newsymbol\curlywedge 1266 
\newsymbol\curlyvee 1267 
\newsymbol\leftthreetimes 1268 
\newsymbol\rightthreetimes 1269 
\newsymbol\subseteqq 136A 
\newsymbol\supseteqq 136B 
\newsymbol\bumpeq 136C 
\newsymbol\Bumpeq 136D 
\newsymbol\lll 136E 
  
\newsymbol\ggg 136F 
  
\newsymbol\circledS 1073 
\newsymbol\pitchfork 1374 
\newsymbol\dotplus 1275 
\newsymbol\backsim 1376 
\newsymbol\backsimeq 1377 
\newsymbol\complement 107B 
\newsymbol\intercal 127C 
\newsymbol\circledcirc 127D 
\newsymbol\circledast 127E 
\newsymbol\circleddash 127F 
\newsymbol\lvertneqq 2300 
\newsymbol\gvertneqq 2301 
\newsymbol\nleq 2302 
\newsymbol\ngeq 2303 
\newsymbol\nless 2304 
\newsymbol\ngtr 2305 
\newsymbol\nprec 2306 
\newsymbol\nsucc 2307 
\newsymbol\lneqq 2308 
\newsymbol\gneqq 2309 
\newsymbol\nleqslant 230A 
\newsymbol\ngeqslant 230B 
\newsymbol\lneq 230C 
\newsymbol\gneq 230D 
\newsymbol\npreceq 230E 
\newsymbol\nsucceq 230F 
\newsymbol\precnsim 2310 
\newsymbol\succnsim 2311 
\newsymbol\lnsim 2312 
\newsymbol\gnsim 2313 
\newsymbol\nleqq 2314 
\newsymbol\ngeqq 2315 
\newsymbol\precneqq 2316 
\newsymbol\succneqq 2317 
\newsymbol\precnapprox 2318 
\newsymbol\succnapprox 2319 
\newsymbol\lnapprox 231A 
\newsymbol\gnapprox 231B 
\newsymbol\nsim 231C 
\newsymbol\ncong 231D 
\newsymbol\diagup 231E 
\newsymbol\diagdown 231F 
\newsymbol\varsubsetneq 2320 
\newsymbol\varsupsetneq 2321 
\newsymbol\nsubseteqq 2322 
\newsymbol\nsupseteqq 2323 
\newsymbol\subsetneqq 2324 
\newsymbol\supsetneqq 2325 
\newsymbol\varsubsetneqq 2326 
\newsymbol\varsupsetneqq 2327 
\newsymbol\subsetneq 2328 
\newsymbol\supsetneq 2329 
\newsymbol\nsubseteq 232A 
\newsymbol\nsupseteq 232B 
\newsymbol\nparallel 232C 
\newsymbol\nmid 232D 
\newsymbol\nshortmid 232E 
\newsymbol\nshortparallel 232F 
\newsymbol\nvdash 2330 
\newsymbol\nVdash 2331 
\newsymbol\nvDash 2332 
\newsymbol\nVDash 2333 
\newsymbol\ntrianglerighteq 2334 
\newsymbol\ntrianglelefteq 2335 
\newsymbol\ntriangleleft 2336 
\newsymbol\ntriangleright 2337 
\newsymbol\nleftarrow 2338 
\newsymbol\nrightarrow 2339 
\newsymbol\nLeftarrow 233A 
\newsymbol\nRightarrow 233B 
\newsymbol\nLeftrightarrow 233C 
\newsymbol\nleftrightarrow 233D 
\newsymbol\divideontimes 223E 
\newsymbol\varnothing 203F 
\newsymbol\nexists 2040 
\newsymbol\Finv 2060 
\newsymbol\Game 2061 
\newsymbol\mho 2066 
\newsymbol\eth 2067 
\newsymbol\eqsim 2368 
\newsymbol\beth 2069 
\newsymbol\gimel 206A 
\newsymbol\daleth 206B 
\newsymbol\lessdot 236C 
\newsymbol\gtrdot 236D 
\newsymbol\ltimes 226E 
\newsymbol\rtimes 226F 
\newsymbol\shortmid 2370 
\newsymbol\shortparallel 2371 
\newsymbol\smallsetminus 2272 
\newsymbol\thicksim 2373 
\newsymbol\thickapprox 2374 
\newsymbol\approxeq 2375 
\newsymbol\succapprox 2376 
\newsymbol\precapprox 2377 
\newsymbol\curvearrowleft 2378 
\newsymbol\curvearrowright 2379 
\newsymbol\digamma 207A 
\newsymbol\varkappa 207B 
\newsymbol\Bbbk 207C 
\newsymbol\hslash 207D 
\undefine\hbar 
\newsymbol\hbar 207E 
\newsymbol\backepsilon 237F 
\catcode`\@=\csname pre amssym.tex at\endcsname 
 
\magnification=1200 
\hsize=468truept 
\vsize=646truept 
\voffset=-10pt 
\parskip=4pt 
\baselineskip=14truept 
\count0=1 
 
\dimen100=\hsize 
 
\def\leftill#1#2#3#4{ 
\medskip 
\line{$ 
\vcenter{ 
\hsize = #1truept \hrule\hbox{\vrule\hbox to  \hsize{\hss \vbox{\vskip#2truept 
\hbox{{\copy100 \the\count105}: #3}\vskip2truept}\hss } 
\vrule}\hrule} 
\dimen110=\dimen100 
\advance\dimen110 by -36truept 
\advance\dimen110 by -#1truept 
\hss \vcenter{\hsize = \dimen110 
\medskip 
\noindent { #4\par\medskip}}$} 
\advance\count105 by 1 
} 
\def\rightill#1#2#3#4{ 
\medskip 
\line{ 
\dimen110=\dimen100 
\advance\dimen110 by -36truept 
\advance\dimen110 by -#1truept 
$\vcenter{\hsize = \dimen110 
\medskip 
\noindent { #4\par\medskip}} 
\hss \vcenter{ 
\hsize = #1truept \hrule\hbox{\vrule\hbox to  \hsize{\hss \vbox{\vskip#2truept 
\hbox{{\copy100 \the\count105}: #3}\vskip2truept}\hss } 
\vrule}\hrule} 
$} 
\advance\count105 by 1 
} 
\def\midill#1#2#3{\medskip 
\line{$\hss 
\vcenter{ 
\hsize = #1truept \hrule\hbox{\vrule\hbox to  \hsize{\hss \vbox{\vskip#2truept 
\hbox{{\copy100 \the\count105}: #3}\vskip2truept}\hss } 
\vrule}\hrule} 
\dimen110=\dimen100 
\advance\dimen110 by -36truept 
\advance\dimen110 by -#1truept 
\hss $} 
\advance\count105 by 1 
} 
\def\insectnum{\copy110\the\count120 
\advance\count120 by 1 
}

\font\ninerm=cmr9 
\font\eightrm=cmr8

\font\tenrm=cmr10 at 10pt 
 
\font\sc=cmcsc10

\def\msb{\fam\msbfam\tenmsb}

\def\bbc{{\msb C}}

\def\bbr{{\msb R}}

\def\bbz{{\msb Z}}

\def\grD{\Delta}

\def\grG{\Gamma}

\def\grS{\Sigma}

\def\gra{\alpha}

\def\gri{\iota} 
 
\def\grl{\lambda} 
 
\def\grn{\nu} 
\def\gro{\omega}

\def\grr{\rho} 
 
\def\grt{\tau}

\def\la#1{\hbox to #1pc{\leftarrowfill}} 
\def\ra#1{\hbox to #1pc{\rightarrowfill}} 
 
\def\fract#1#2{\raise4pt\hbox{$ #1 \atop #2 $}} 
\def\decdnar#1{\phantom{\hbox{$\scriptstyle{#1}$}} 
\left\downarrow\vbox{\vskip15pt\hbox{$\scriptstyle{#1}$}}\right.} 
\def\decupar#1{\phantom{\hbox{$\scriptstyle{#1}$}} 
\left\uparrow\vbox{\vskip15pt\hbox{$\scriptstyle{#1}$}}\right.} 
\def\bowtie{\hbox to 1pt{\hss}\raise.66pt\hbox{$\scriptstyle{>}$} 
\kern-4.9pt\triangleleft} 
\def\hsmash{\triangleright\kern-4.4pt\raise.66pt\hbox{$\scriptstyle{<}$}} 
\def\boxit#1{\vbox{\hrule\hbox{\vrule\kern3pt 
\vbox{\kern3pt#1\kern3pt}\kern3pt\vrule}\hrule}}

\def\za{\vrule height6pt width4pt depth1pt}

\font\aa=eufm10

\def\Got#1{\hbox{\aa#1}}

\def\bfa{{\bf a}}

\def\cald{{\cal D}} 
 
\def\calf{{\cal F}}

\def\call{{\cal L}}

\def\calz{{\cal Z}}

\def\gc{{\Got c}}

\def\gg{{\Got g}} 
\def\gh{{\Got h}}

\def\gs{{\Got s}} 
\def\gt{{\Got t}}

\def\gC{{\Got C}}

\def\gG{{\Got G}} 
\def\gH{{\Got H}}

\def\gS{{\Got S}} 
\def\gT{{\Got T}}

\def\Got#1{\hbox{\aa#1}}

\def\heta{\hat{\eta}} 
\def\gsp1{{\Got s}{\Got p}(1)}

\font\svtnrm=cmr17

\font\bsc=cmcsc10 at 10truept

\def\teta{\tilde{\eta}}
\def\tmu{\tilde{\mu}}
\def\hgt{\hat{\gt}}

\def\hgrD{\hat{\grD}}
\def\hgrl{\hat{\grl}}
\def\hgra{\hat{\gra}}
\def\hf{\hat{f}}
\def\hell{\hat{\ell}}
\def\heta{\hat{\eta}}
\def\teta{\tilde{\eta}}
\def\tM{\tilde{M}}
\def\hgro{\hat{\gro}}
\def\hphi{\hat{\phi}}
\def\hpi{\hat{\pi}}

\centerline{\svtnrm A Note on Toric Contact Geometry}
\bigskip
\centerline{\sc Charles P. Boyer~~ Krzysztof Galicki~~}
\footnote{}{\ninerm During the preparation of this work the authors 
were partially supported by NSF grant DMS-9970904.}
\bigskip
\centerline{\vbox{\hsize = 5.85truein
\baselineskip = 12.5truept
\eightrm
\noindent {\bsc Abstract:}
After observing that the well-known convexity theorems of symplectic geometry
also hold for compact contact manifolds with an effective torus action whose
Reeb vector field corresponds to an element of the Lie algebra of the torus,
we use this fact together with a recent symplectic orbifold
version of Delzant's theorem due to Lerman and Tolman 
[LT] to show that every such compact toric contact manifold can be
obtained by a contact reduction from an odd dimensional sphere.}}
\tenrm

\bigskip
\baselineskip = 10 truept
\centerline{\bf 1. Introduction} 
\bigskip

The main purpose of this note is to prove a Delzant-type
Theorem [Del] which says that every compact toric contact manifold whose
Reeb vector field corresponds to an element of the Lie algebra of the torus
(a condition we call {\it of Reeb type}) can be obtained by contact reduction
from an odd dimensional sphere. This result makes use of the symplectic
orbifold version of Delzant's Theorem by Lerman and Tolman [LT] thereby
showing the close relationship between the two. Toric contact geometry has
been previously considered by Banyaga and Molino [BM1,BM2]. They show that
there are two cases: (1) The action of the torus is regular in which case the
image of the moment map is a sphere. Thus, for $n\geq 3$ the original contact
manifold is the product $T^{n+1}\times S^n.$ (2) The action is singular in
which case the image generates a closed convex polytope. Then Banyaga and
Molino give a Delzant-type theorem by showing that this polytope determines
the toric contact structure up to isomorphism. But they do not prove that
every such manifold can be obtained from reduction. Indeed, this may not be
true. However, if one makes the added assumption that the Reeb vector field
corresponds to an element of the Lie algebra of the torus, such a result is
true as we shall show. While [BM1,BM2] do consider this case they do not prove
such a result since their proofs are entirely different from the usual proofs
of the Delzant theorems in symplectic geometry [Del, Gui, LT] which make use
of reduction.

We begin by reviewing some well known facts about contact manifolds and
symplectic cones. We then discuss the moment maps associated to toral actions. 
In particular, we show that the well-known convexity theorem of Atiyah [At],
Guillemin and Sternberg [GS] hold for compact contact manifolds of Reeb
type, a result given previously in the toric case by Banyaga and Molino [BM1,BM2].
In contact geometry we need to fix a contact 1-form, that is a Pfaffian
structure, within the contact structure. Doing so we show that a torus that
acts effectively on a compact contact manifold and preserves the Pfaffian form
gives rise to a moment map whose image is a convex polytope lying in a certain
hyperplane, which we call the {\it characteristic hyperplane}, in the dual of
the Lie algebra of the torus. If one changes the Pfaffian form by a function
that is invariant under the torus, the resulting polytope then differs from
the former by a change of scale. Indeed, it is always possible to choose the
Pfaffian form so that the polytope is rational. Following Lerman and Tolman we
also give a theorem  relating the geometry of labeled polytopes to the
geometry of the contact manifold with a fixed Pfaffian structure with its torus
action and characteristic foliation.

Our motivation for this study is that of obtaining explicit positive Einstein
metrics using methods of contact geometry that the authors have recently
developed [BGM, BGMR, BG1-2]. In a forthcoming work we investigate
precisely which toric contact manifolds admit Sasakian-Einstein metrics.

\noindent{\sc Acknowledgments}: We would like to thank E. Lerman for pointing
our some incorrect statements in the first draft of this paper. The second
author would like to thank Max-Planck-Institute f\"ur Mathematik in Bonn for
support and hospitality. This article was completed during his two month stay
in Bonn in 1999.

\bigskip
\centerline {\bf 2. Symplectic Cones and Contact Manifolds}
\bigskip

\noindent{\sc Definition} 2.1 \tensl A {\it contact structure} on a manifold
$M$ is a subbundle $\cald,$ called the contact distribution, of the tangent
bundle $TM$ that is maximally non-integrable. \tenrm

When $M$ is orientable, which we shall always assume, the contact distribution
can be given as the kernel of a 1-form $\eta$ that satisfies the nondegeneracy
condition
$$ \eta\wedge (d\eta)^n\neq 0 \leqno{2.2}$$
and $d\eta$ is a symplectic form on the subbundle $\cald.$
It should be realized, however, that the contact structure is defined not by
$\eta$, but by an equivalence class of such 1-forms. Explicitly, two such
1-forms $\eta,\eta'$ are equivalent if there exists a nowhere vanishing
function $f$ on $M$ such that $\eta'=f\eta.$ For us it will often be convenient
to fix a contact form $\eta$ in the equivalence class. In [LM] this is referred
to as a {\it Pfaffian structure}.

\noindent{\sc Definition} 2.3: \tensl A {\it symplectic cone} is a cone
$C(M)=M\times \bbr^+$ with a symplectic structure $\gro$ with a one parameter
group $\grr_t$ of homotheties of $\gro$ whose infinitesimal generator is a
vector field on $\bbr^+.$ 
\tenrm

Definition 2.3 is referred to as a symplectic Liouville structure in [LM].
We denote by $r$ the coordinate on $\bbr^+$ in which case the infinitesimal
generator of the homothety group is the Liouville vector field $\Psi =
r\partial_r.$ Define the one form $\teta$ on $C(M)$ by $\teta =\Psi \rfloor
\gro.$ As $\Psi$ is an infinitesimal
generator of homotheties of $\gro$ 
it follows that $\gro$ is the exact 2-form
$\gro=d\teta.$ Let us identify $M$ with $M\times \{1\},$ and define $\eta
=\teta|_{M\times \{1\}}.$ Then we see that $\teta =r\eta,$ and this gives
$\gro$ as $$\gro=dr\wedge \eta +rd\eta.  \leqno{2.4}$$ The condition that
$\gro^{n+1}\neq 0$ implies that $\eta$ gives $M$ the structure of a contact
manifold with a fixed 1-form, {\it i.e.,}  2.2 is satisfied.  Conversely, one can
easily reverse this procedure: given a contact manifold with a fixed 1-form
$(M,\eta),$ defining $\gro$ on $C(M)$ by 2.4 gives the cone $C(M)$ a
symplectic structure with homotheties.  In [Eli, MS] $(C(M),\gro)$ 
is called the {\it
symplectization} of $M$ while it is called
the {\it symplectification} of $M$ in [LM].
The terminology is apparently due to Arnold.  We have arrived at the well-known
result

\noindent{\sc Proposition} 2.5: \tensl An orientable manifold (orbifold) is
a contact manifold (orbifold) if and only if the cone $C(M)=M\times \bbr^+$ is
a symplectic cone. 
\tenrm

Now a contact manifold with a fixed 1-form $\eta$ has associated to it a unique
vector field called the {\it characteristic} or Reeb vector field $\xi$ defined
by the conditions
$$\eta(\xi)=1,\qquad \hbox{and}\qquad \xi\rfloor d\eta=0.\leqno{2.6}$$ 
Clearly $\xi$ is nowhere vanishing on $M$ and thus induces a foliation
$\calf^\xi$ of $M$ called the {\it characteristic foliation}. It should be
mentioned that both $\xi$ and its characteristic foliation depend on the choice
of 1-form $\eta,$ that is, they are invariants of the Pfaffian structure, but
not the contact structure. Conversely, a characteristic vector field uniquely
determines the 1-form $\eta$ within the contact structure $\cald.$

\bigskip
\centerline{\bf 3. Contact Transformations}
\bigskip

Let $M$ be an orientable contact manifold, and let $\gC(M,\cald)$ denote the
group of contact transformations, {\it i.e.,} the subgroup of the group
$\hbox{Diff}(M)$ of diffeomophisms of $M$ that leaves the contact distribution
$\cald$ invariant. If we fix a contact form $\eta$ such that
$\cald=\hbox{ker}~\eta,$ then $\gC(M,\cald)$ can be characterized as the
subgroup of diffeomorphisms $\phi:M\ra{1.3} M$ that satisfy $\phi^*\eta =f\eta$
for some nowhere vanishing $f\in C^\infty(M).$ With the 1-form $\eta$ fixed we
are interested in the subgroup $\gC(M,\eta)$ of {\it strict contact
transformations} defined by the condition $\phi^*\eta =\eta.$ Likewise, we
denote the ``Lie algebras'' of $\gC(M,\cald)$ and $\gC(M,\eta)$ by
$\gc(M,\cald)$ and $\gc(M,\eta),$ respectively. They can be characterized as
follows:
$$\gc(M,\cald)=\{X\in \grG(M)|\call_X\eta=g\eta~ \hbox{for some}~ g\in
C^\infty(M)\},$$
$$\gc(M,\eta)=\{X\in \grG(M)|\call_X\eta =0\},$$
where $\grG(M)$ denotes the Lie algebra of smooth vector fields on $M.$ Clearly,
$\gc(M,\eta)$ is a Lie subalgebra of $\gc(M,\cald).$ In contrast to the
symplectic case, every infinitesimal contact transformation is Hamiltonian.
More explicitly,

\noindent{\sc Proposition} 3.1: [LM] \tensl The contact 1-form $\eta$ induces
a Lie algebra isomorphism $\Phi$ between the Lie algebra $\gc(M,\cald)$ of
infinitesimal contact transformations and the Lie algebra $C^\infty(M)$ of
smooth functions with the Jacobi bracket defined by $$\Phi(X)=\eta(X).$$
Furthermore, under this isomorphism the subalgebra $\gc(M,\eta)$ is isomorphic
to the subalgebra $C^\infty(M)^\xi$ of functions in $C^\infty(M)$ that are
invariant under the flow generated by the Reeb vector field $\xi.$ In
particular, $\xi =\Phi^{-1}(1).$ \tenrm

\noindent The function $\eta(X)$ is known as a {\it contact Hamiltonian
function}.

Similarly, on $C(M)$ we let $\gS(C(M),\gro)$ denote the group of
symplectomorphisms of $(C(M),\gro),$ and $\gS_0(C(M),\gro)$ the subgroup of
$\gS(C(M),\gro)$ that commutes with homotheties, 
{\it i.e.,} the automorphism group of
the symplectic Liouville structure. The corresponding Lie algebras are denoted
by $\gs(C(M),\gro)$ and $\gs_0(C(M),\gro),$ respectively. Given a vector field
$X\in \gs_0(C(M),\gro),$ the fact that $X$ commutes with the Liouville
vector field $\Psi$ implies that $X$ is projectable to a vector field $X_M$ on
$M$ and one easily sees that [LM]

\noindent{\sc Proposition} 3.2: \tensl The maps $X\mapsto X_M\mapsto
\eta(X_M)$ induce Lie algebra isomorphisms $\gs_0(C(M),\gro)\approx
\gc(M,\eta)\approx C^\infty(M)^\xi.$ Furthermore, $\xi$ is in the center
of $\gc(M,\eta).$ \tenrm

\bigskip
\centerline{\bf 4. Convexity and the Moment Map}
\bigskip

We now consider the moment map construction for $C(M)$ and $M.$ Let $\gG$
be a Lie group acting on the symplectic cone $(C(M),\gro)$
which leaves invariant the symplectic form $\gro$ and which commutes with the
homothety group. In particular, $\gG$ is a subgroup of the group
$\gS_0(C(M),\gro)$ and gives rise to a moment map
$$\tmu:C(M)\ra{1.5} \gg^*, \leqno{4.1}$$
where $\gg$ denotes the Lie algebra of $\gG$ and $\gg^*$ its
dual. Explicitly $\tmu$ is defined by 
$$d<\tmu,\grt>= -X^\grt\rfloor \gro, \leqno{4.2}$$
where $X^\grt$ denotes the vector field on $C(M)$ corresponding to $\grt\in
\gg.$ For simplicity we denote the function $<\tmu,\grt>$ by $\tmu^\grt.$
An easy computation then shows that up to an additive constant
$$\tmu^\grt=\teta(X^\grt). \leqno{4.3}$$ 
With this choice $\tmu$ is clearly homogeneous. Indeed, under homotheties
$r\mapsto e^tr$ we have
$$\gro\mapsto e^t\gro, \qquad \teta\mapsto e^t\teta, \qquad \tmu\mapsto
e^t\tmu.$$
Now consider the case when the Lie group $\gG$ is a torus $\gT$ with Lie
algebra $\gt.$ Then there is the following convexity theorem [deMT]:

\noindent{\sc Theorem} 4.4~[deMT]: \tensl Let $(C(M),\gro)$ be a symplectic
cone with $M$ compact, and let $\gT\subset \gS_0(C(M),\gro)$ be a torus group.
Assume further that there is an element $\grt\in \gt$ such that $<\mu,\grt>
>0.$ Then the image $\mu(C(M))$ is a convex polyhedral cone.  \tenrm

Now let us return to the case of an arbitrary Lie group $\gG\subset
\gS_0(C(M),\gro).$ Again we identify $M$ with $M\times \{1\}.$ Since $\gG$
commutes with homotheties we get an induced action of $\gG$ on $M,$ and by
Proposition 3.2 we can identify $\gG$ with a Lie subgroup of $\gC(M,\eta).$
This gives a moment map by restriction, viz.
$$\mu:M\ra{1.5} \gg^*,$$
$$\mu=\tmu|_{M\times \{1\}}, \qquad \mu^\grt=\eta(X^\grt). \leqno{4.5}$$
Such a moment map was noticed previously by Geiges [Gei] and in [BGM] within
the context of 3-Sasakian geometry.

We wish to consider the special case when the Lie group $\gG$ is an
$n+1$-dimensional torus $\gT^{n+1}.$ Let $\gt_{n+1}$ denote the Lie algebra of
$\gT^{n+1},$ and let $\{e_i\}_{i=0}^n$ denote the standard basis for
$\gt_{n+1}\approx \bbr^{n+1}.$ Corresponding to each basis element $e_i$ there
is a vector field $X^{e_i}$ which for convenience we denote by $H_i.$  We shall
assume that the Reeb vector field $\xi$ corresponds to an element $\varsigma,$
which we call the {\it characteristic vector}, in the Lie algebra $\gt_{n+1}.$
Hence, the Reeb vector field is almost periodic. In this case we say that the
torus action is {\it of Reeb type}. Let $\{e_i^*\}_{i=0}^n$ denote the dual
basis of $\gt_{n+1}^*.$ Then we can write the corresponding characteristic
vector $\varsigma$ and moment map $\mu$ as 
$$\varsigma =\sum_{i=0}^na_ie_i, \qquad \mu= \sum_{i=0}^n\mu_ie_i^*
\leqno{4.6}$$ 
for some $a_i\in \bbr.$ Then the moment map $\mu$ maps
$\varsigma$ to the hyperplane given by 
$$<\mu,\varsigma> =\sum_{i=0}^n\mu_ia_i
=1.  \leqno{4.7}$$ 
We call this hyperplane {\it the characteristic hyperplane}.  The nondegeneracy
of $\eta$ implies that the plane defined by 4.7 is actually a hyperplane,
{\it i.e.,}  its codimension is one. Now there is a commutative diagram
$$\matrix{C(M)&&\cr
          \decupar{}&\fract{\tmu}{\searrow}&\cr
          M&\fract{\mu}{\ra{1.8}}&\gt_{n+1}^*,} \leqno{4.8}$$
where the vertical arrow is the natural inclusion $M$ into $C(M)$ as $M\times
\{1\}.$ Thus, intersecting the characteristic hyperplane 4.6 with the convex
cone of Theorem 4.4 we arrive at:

\noindent{\sc Theorem} 4.9: \tensl Let $(M,\eta)$ be a compact contact
manifold with a fixed contact form $\eta,$ and let $\mu$ be the moment map of a
torus $\gT$ acting effectively on $M$ which preserves the contact form $\eta.$
Suppose also that the Reeb vector field $\xi$ corresponds to an element of the
Lie algebra $\gt$ of $\gT,$ so that the torus action is of Reeb type. Then the
image $\mu(M)$ is a convex compact polytope lying in the characteristic
hyperplane.  \tenrm

This theorem was proved earlier by a different method by Banyaga and Molino
[BM1,BM2], and they also give an example due to R. Lutz which shows that the
assumption on the Reeb vector field is necessary.  

We are mainly interested in rational polytopes. Let
$\ell\subset \gt$ denote the lattice of circle subgroups of $\gT.$ We recall
some definitions for convex polytopes [Zie, Gui, LT]

\noindent{\sc Definition} 4.10: \tensl A {\it facet} is a codimension
one face. A convex polytope of dimension $n$ is called {\it simple} if there
are precisely $n$ facets meeting at each vertex. A convex polytope $\grD\subset
\gt^*$ is {\it rational} if, for some $\grl_i\in \bbr$, there are $y_i\in \ell$
such that 
$$\grD= \bigcap_{i=1}^N\{\gra\in \gt^*~|~ <\gra,y_i>\leq \grl_i\},$$
where $N$ is the number of facets of $\grD.$ \tenrm

We are now ready for

\noindent{\sc Theorem} 4.11: \tensl Under the hypothesis of Theorem 4.9,
the polytope $\mu(M)\subset \gt^*$ is simple of dimension $\dim \gt -1,$ and
it is rational if and only if the characteristic vector $\varsigma$ lies in the
lattice $\ell\subset \gt$ of circle subgroups of $\gT.$ \tenrm

\noindent{\sc Proof}: The nondegeneracy of $\eta$ implies that the
dimension of the polytope $\grD$ is $\dim~\gt-1,$ and the fact that $\grD$ is
the intersection of a convex polyhedral cone and a hyperplane implies that
$\grD$ is simple. To verify the rationality condition consider the characteristic
foliation $\calf^\xi.$ The space of leaves $\calz$ is a compact orbifold if and
only if $\varsigma$ is in the lattice $\ell.$ So by the quasi-regular version
[Tho] of the Boothby-Wang fibration, this orbifold has a symplectic structure.
Moreover, since $\gt$ is Abelian and contains $\varsigma,$ there is a Lie
algebra $\hat{\gt}$ acting on $\calz$ by infinitesimal symplectomorphisms which
fits into the exact sequence 
$$0\ra{1.5}\{\varsigma\}\ra{1.5}\gt\fract{\varpi}{\ra{1.5}}\hat{\gt}\ra{1.5}
0,$$
where $\{\varsigma\}$ denotes the Lie algebra generated by $\varsigma.$
Consider the dual sequence
$$0\ra{1.5}\hat{\gt}^*\fract{\varpi^*}{\ra{1.5}}\gt^*
{\ra{1.5}}\{\varsigma\}^* \ra{1.5} 0.\leqno{4.12}$$
By a theorem of Lerman and Tolman [LT] there is a moment map $\hat{\mu}:\calz
\ra{1.3} \hgt^*$ whose image in $\hgt^*$ is a rational convex polytope
$\hat{\grD}.$  Furthermore, if $\pi:M\ra{1.3} \calz$ denotes the orbifold
V-bundle map and $\gro$ the symplectic form on $\calz$ we have that
$\pi^*\gro=d\eta.$ It follows that $\varpi^*\pi^*\hat{\mu}$ differs from $\mu$
by a constant $c\in \gt^*,$ that is
$$\mu=\varpi^*\pi^*\hat{\mu} +c. \leqno{4.13}$$
Moreover, since $\mu$ maps into the characteristic
hyperplane $c$ satisfies $<c,\varsigma>=1.$ Let $\grD$ denote the polytope
$\mu(M).$ From Definition 4.10 we see that the polytope $\grD$ is rational if
and only if we can choose the $y_i$ to lie in $\ell.$ But since $\hgrD$ is
rational and $\gra\in \grD$ if and only if $\gra -c\in \hgrD,$ there are
$y_i\in \ell$ such that $\hgrD$ is determined by $<\hgra,y_i>\leq \hgrl_i.$
Thus, $\grD$ is given by the equation
$$\bigcap_{i=1}^N\{\gra\in \gt^*~|~<\gra,y_i>~\leq~ \grl_i=\hgrl_i+<c,y_i>\}$$
which is clearly rational. \hfill\za

For quasi-regular contact Pfaffian structures equation 4.13 gives the
fundamental relation between our contact moment map and the moment map for
compact symplectic orbifolds described by Lerman and Tolman. Moreover, in [LT]
it is shown that the polytope $\hgrD$ has labels associated with each facet.
The outward normal vector $y_i$ to the $i^{\rm th}$ facet $\hf_i$ of $\hgrD$ lies in
$\ell,$ so there is a positive integer $m_i$ and a primitive vector $p_i\in
\ell$ such that $y_i=m_ip_i.$ Thus, by associating $m_i$ to the $i^{\rm th}$
facet $f_i$
for each $i=1,\ldots, N$ we obtain $\hgrD$ as a labeled polytope. But from the
discussion above the outward normals $y_i$ to the $i^{\rm th}$ facet $f_i$ of $\grD$
coincide with the outward normals to $\hf_i$ of $\hgrD,$ so that $\grD$ is also
a labeled polytope with the integers $m_i$ associated to the 
$i^{\rm th}$ facet $f_i.$
Thus an immediate consequence of 4.13 and Lemma 6.6 of [LT] is

\noindent{\sc Theorem} 4.14: \tensl Let $(M,\eta)$ be a contact manifold
with a fixed quasi-regular contact form $\eta$ with an effective action of a
torus $\gT$ of Reeb type that leaves the 1-form $\eta$ invariant,  For every
point $x\in M$ let $F(x)$ denote the set of open facets of $\grD$ whose closure
contains $\mu(x)$ and let $m_i$ and $p_i$ denote the labels and primitive
outward normal vectors to the $i^{\rm th}$ facet.  Then the Lie algebra
$\gh_x$ of the isotropy subgroup $\gH_x$ of $\gT$ at $x$ is the linear span of
the vectors $p_i$ for all $i$ such that the $i^{\rm th}$
open facet $f_i$ lies in $F(x).$ In particular, if $\mu(x)$
is a vertex of the polytope $\grD,$ then $\gH_x$ is isomorphic to the factor
group $\gT/S^1(\xi)$ where $S^1(\xi)$ denotes the circle subgroup generated by
the characteristic vector field $\xi.$ Furthermore, $\grD$ is the convex hull
of its vertices.

Let $\ell_x\in \gh_x$ denote the lattice of circle subgroups of $\gH_x$ and let
$\hell_x$ denote the sub-lattice of $\ell$ generated by the vectors
$\{m_ip_i\}_{f_i\in F(x)}.$ Then the leaf holonomy group at $x$ of the
characteristic foliation $\calf^\xi$ is isomorphic to $\ell_x/\hell_x.$ In
particular, $(M,\eta)$ is regular if and only if the set $\{m_ip_i\}_{f_i\in
F(x)}$ generates $\ell_x$ for all $x\in M$. (In this case $m_i=1\ \ \forall i$).
\tenrm

We now consider varying the contact form $\eta$ within the contact structure.
Let $\eta'=f\eta$ where $f$ is a nowhere vanishing function on $M,$ and let
$\xi'$ be the Reeb vector field associated to $\eta'.$ Write $\xi'=\xi +\rho.$
Suppose further that $X\in \gc(M,\eta).$ Then we have the following elementary 
lemma whose proof we leave to the reader:

\noindent{\sc Lemma} 4.15: \tensl The following hold:
\item{(i)} $\displaystyle{f= {1\over \eta(\grr)+1}}.$
\item{(ii)}$\grr\rfloor d\eta= -d(\eta(\grr))+\xi(\eta(\grr))\eta.$
\item{(iii)} $\displaystyle{\call_\grr\eta =\call_{\xi'}\eta
=\xi(\eta(\grr))\eta},$ {\it i.e.,} 
both $\xi'$ and $\grr$ are infinitesimal contact
transformations.
\item{(iv)} $\xi',\grr\in \gc(M,\eta)$ if and only if $f\in C^\infty(M)^\xi.$
\item{(v)} $X\in \gc(M,\eta')$ if and only if $Xf=0.$
\tenrm

We return to the case when $(M,\eta)$ is a contact manifold with an
$(n+1)$-torus $\gT^{n+1}$ acting as strict contact transformations whose
characteristic vector $\varsigma$ lies in the Lie algebra $\gt_{n+1}.$  The
moment map for this torus action is given by 4.6 with $\mu_i=\eta(H_i).$ If
$\eta'=f\eta$ is another 1-form in same contact structure where $f$ is
invariant under the torus action, then this action preserves $\eta'$ as well,
and the corresponding moment map $\mu'$ satisfies $\mu'=f\mu.$ Let $\grD$ be
the polytope associated with $\eta.$  Then the polytope $\grD'$ associated to
$\eta'$ has the same dimension as $\grD$ with the same number of facets and the
same number of vertices. The size of the faces and labels, however, depend on
the contact form $\eta,$ that is on the Pfaffian structure, and the labels are
defined only when $\eta$ is quasi-regular. But vectors lying in the lattice
$\ell$ of circle subgroups are dense in $\gt_{n+1},$ so from the point of view
of the contact structure we can always perturb the Reeb vector field and
contact form so that the characteristic foliation is quasi-regular, and by
Theorem 4.11 so that the polytope is rational. We are ready for

\noindent{\sc Definition} 4.16: \tensl A contact manifold (orbifold)
$(M,\cald)$ of dimension $2n+1$ is called a {\it toric contact manifold
(orbifold)} (written as the triple $(M,\cald,\gT)$) if there are a 1-form
$\eta$ that represents the contact structure $\cald,$ and an effective action
of an $(n+1)$-dimensional torus $\gT$ on $M$ that preserves the contact form
$\eta.$  If in addition the Reeb vector field associated to $\eta$ corresponds
to an element of the Lie algebra $\gt$ of $\gT,$ we say that $(M,\cald,\gT)$
is a {\it toric contact manifold of Reeb type}. \tenrm
 
If we wish to fix a contact form $\eta$ we write $(M,\eta,\gT)$ for a toric
contact manifold instead of $(M,\cald,\gT).$  Our discussion above proves

\noindent{\sc Proposition} 4.17: \tensl Let $(M,\cald,\gT)$ be a compact toric
contact orbifold of Reeb type, then $\cald$ can be represented by a
quasi-regular contact form, and hence, by a rational polytope. \tenrm

The fiduciary examples of compact contact manifolds are the odd dimensional
spheres $S^{2n+1}.$ 

\noindent{\sc Example} 4.18: $S^{2n+1}$ with the standard contact
structure. This is the contact structure induced from the standard symplectic
structure on $\bbr^{2n+2}$ given in Cartesian coordinates
$(x_0,y_0,\ldots,x_n,y_n)\in \bbr^{2n+2}$ by 
$$\gro=2\sum_{i=0}^ndx_i\wedge dy_i.$$
In this case the standard contact form $\eta$ and the standard characteristic
vector field $\xi$ are given by
$$\eta = {1\over r}\sum_{i=0}^n(x_idy_i-y_idx_i),\qquad
\xi=\sum_{i=0}^n(x_i\partial_{y_i}-y_i\partial_{x_i}),$$
where $r=\sum_{i=0}^n(x_i^2+y_i^2)$ (not the usual $r$). The maximal torus
$\gT^{n+1}$ is generated by the vector fields $H_i=
x_i\partial_{y_i}-y_i\partial_{x_i}$ for $i=0,\ldots, n,$ so $(S^{2n+1},\eta)$
is a toric contact manifold. The moment map is easily seen to be
$$\mu(x_0,y_0,\ldots,x_n,y_n)= {1\over r}\sum_{i=0}^n(x_i^2+y_i^2)e_i^*.$$
Letting $r_0,\ldots,r_n$ denote the coordinates for $\gt_{n+1}^*$ we see that
the characteristic hyperplane is just
$$r_0+\cdots +r_n=1.$$
Thus, the image $\mu(S^{2n+1})$ is just the standard $n$-simplex with $r_i$'s
as barycentric coordinates.

Now we can consider non-standard characteristic vector fields and Pfaffian
forms within the standard contact structure on $S^{2n+1}.$ These are
deformations [YK] of the standard form depending on $n+1$ positive real
parameters $(a_0,\ldots,a_n)\in (\bbr^+)^{n+1}.$ In this case the 1-form
$\eta$ and characteristic vector fields are given by
$$\eta_\bfa = {\sum_{i=0}^n(x_idy_i-y_idx_i)\over\sum_{i=0}^n
a_i(x_i^2+y_i^2)} \qquad \xi_\bfa
=\sum_{i=0}^na_i(x_i\partial_{y_i}-y_i\partial_{x_i}),$$
so that 
$$\eta_\bfa =\biggl({r\over \sum_{i=0}^n a_i(x_i^2+y_i^2)}\biggr)\eta.$$
The characteristic hyperplane is $$\sum_{i=0}^na_ir_i=1,$$ so the polytope is
given by the ``weighted''$n$-simplex determined by this and $$0\leq r_i\leq
{1\over a_i}.$$ The special case where $a_i=a$ for all $i=0,\ldots,n$ is just
the dilated standard $n$-simplex 
$$a(r_0+\cdots+r_n)=1, \qquad 0\leq r_i\leq {1\over a}.$$

\bigskip
\centerline{\bf 5. A Delzant Theorem for Toric Contact Manifolds of Reeb type}
\bigskip

We begin by considering contact reduction [BGM, Gei]. Let $(\tM,\teta)$ be a
compact contact manifold with a fixed quasi-regular contact form $\teta.$
Suppose also that a compact Lie group $\gG$ acts on $\tM$ preserving the
contact form $\teta$ and let $\mu:\tM\ra{1.3} \gg^*$ denote the corresponding
moment map. Then if $\gG$ acts freely on the zero set
$\mu^{-1}(0)\fract{\gri}{\hookrightarrow} \tM,$ the quotient
$M=\mu^{-1}(0)/\gG$ is a compact contact manifold with a unique fixed 1-form
$\eta$ satisfying $\gri^*\teta =p^*\eta$ where $p:\mu^{-1}(0)\ra{1.3} M$
denotes the natural projection.

We shall prove

\noindent{\sc Theorem} 5.1: \tensl Let $(M,\eta,\gT)$ be a compact toric
contact manifold of Reeb type and a fixed quasi-regular
contact form $\eta.$ Then $(M,\eta)$ is isomorphic to the reduction by a torus
of a sphere $S^{2N-1}$ with its standard contact structure and with a fixed
1-form $\eta_\bfa$. \tenrm  

\noindent{\sc Proof}: Now $M$ is compact of dimension, say, $2n+1,$ and 
since $\eta$ is quasi-regular, the space of leaves $\calz$ of the
characteristic foliation $\calf$ is a compact symplectic orbifold of dimension
$2n.$ Furthermore, since $M$ is toric, so is $\calz,$ that is, there is an $n$
dimensional torus $\gT^n$ preserving the symplectic structure $\gro$ on
$\calz.$ Furthermore, by [BG2] $\gro$ represents an integral class in
$H^2_{orb}(\calz,\bbz)$ and $M$ is the total space of the principal $S^1$
V-bundle $M\fract{\pi}{\ra{1.3}} \calz$ whose first Chern class is represented
by $d\eta =\pi^*\gro.$

Now by Theorem 8.1 of Lerman and Tolman [LT] $(\calz,\gro)$ is isomorphic to
the symplectic reduction $(\bbc^N,\gro_0)$ with the standard symplectic
structure by a torus $\gT^{N-n}$of dimension $N-n.$ If
$\mu_{N-n}:\bbc^N\ra{1.8} \gt^*_{N-n}$
denotes the moment map for the $\gT^{N-n}$ action, then $(\calz,\gro)$ is
isomorphic to $(\mu_{N-n}^{-1}(\grl)/\gT^{N-n},\hgro)$ where $\grl$ is a
regular value of $\mu_{N-n}$ and $\hgro$ is the unique symplectic 2-form
induced by reduction. Let $\hphi$ denote the above isomorphism. It follows
that the cohomology class of $\hgro =(\hphi^{-1})^*\gro$ is integral in
$H^2_{orb}(\mu_{N-n}^{-1}(\grl)/\gT^{N-n},\bbz).$ By the orbifold
version of the Boothby-Wang theorem there is an $S^1$ V-bundle  $\pi:P\ra{1.3}
\mu_{N-n}^{-1}(\grl)/\gT^{N-n},$ a connection form $\heta$ on $P$ such that
$d\heta = \hpi^*\hgro.$  Thus there is an $S^1$-equivariant V-bundle map 
$$\matrix{M&\fract{\phi}{\ra{3.8}}& P&\cr
            &&\cr
           \decdnar{\pi} && \decdnar{\hpi} &\cr
            &&\cr
            \calz &\fract{\hphi}{\ra{1.8}}&
\mu_{N-n}^{-1}(\grl)/\gT^{N-n}& \cr} \leqno{5.2}$$
such that $\phi^*d\heta=d\eta.$ Thus, 
$\phi^*\heta$ and $\eta$ differ by a closed 1-form. But the space of closed
1-forms on $M$ is path-connected, so 
one can find a one parameter family of
connections having the same curvature that connect $\phi^*\heta$  to $\eta.$ 
So by Gray's stability theorem [MS] $\phi^*\heta$ and $\eta$ define the same
contact structure. Thus, we can choose $\phi^*\heta =\eta.$  Moreover, by
equivariance the characteristic vector $\varsigma$ of the contact manifold
$(P,\heta)$ lies in the Lie algebra $\gt_{N-n},$ so we can split off the
circle that it generates and write $\gT^{N-n}=S^1_\varsigma \times
\gT^{N-n-1}$, where $\gT^{N-n-1}$ is an  $(N-n-1)$-dimensional torus. It
follows that $P= \mu_{N-n}^{-1}(\grl)/\gT^{N-n-1}.$ Hereafter, we identify
$(M,\eta)$ with $(\mu_{N-n}^{-1}(\grl)/\gT^{N-n-1},\heta).$ Now
$\mu_{N-n}^{-1}(\grl)$ is a torus bundle over a compact manifold, so it is a
compact manifold which by construction is an intersection of $N-n$ real
quadrics in $\bbc^N.$ It follows that there is a component of the moment map
$\mu_{N-n}$ which takes the form $\sum_ia_i|z_i|^2$ with $a_i>0$ for all $i.$
Let $\bfa$ denote the vector in $\bbr^N$ whose $i^{\rm th}$ component is
$a_i,$ and consider the ellipsoid $\grS_\bfa=\{\sum_ia_i|z_i|^2=1\}\cong
S^{2N-1}.$ Then there is a $\gT^{N-n}$-moment map $\grn_\bfa:\grS_\bfa\ra{1.5}
t_{N-n}^*$ such that $\grn_\bfa^{-1}(0)=  \mu_{N-n}^{-1}(\grl).$ Furthermore,
letting $\eta_0=\sum_i(x_idy_i-y_idx_i)$ we see that $d\eta_0=\gro_0$ on
$\bbc^N$ and that $\eta_0|_{\grS_\bfa}= \eta_\bfa|_{\grS_\bfa},$ where 
$\eta_\bfa$ is the deformed 1-form of Example 4.18. Thus, letting
$p:\mu_{N-n}^{-1}(\grl)\ra{1.5} M$ denote the natural submersion, and 
$\gri:\mu_{N-n}^{-1}(\grl)\ra{1.5} \grS_\bfa$ the natural inclusion, we see
that $p^*\eta =\gri^*\eta_\bfa,$ so $(M,\eta)$ is obtained from
$(\grS_\bfa,\eta_\bfa)$ by contact reduction. 
\hfill\za

Now Theorem 9.1 of [LT] says that every symplectic toric orbifold possesses
an invariant complex structure which is compatible with its symplectic form.
This means that every symplectic toric orbifold is actually K\"ahler and we can
combine this fact with our results to get the following:

\noindent{\sc Theorem} 5.3: \tensl Every compact toric contact manifold of
Reeb type admits a compatible invariant Sasakian
structure. \tenrm

\bigskip
\bigskip
\medskip
\centerline{\bf Bibliography}
\medskip
\font\ninesl=cmsl9
\font\bsc=cmcsc10 at 10truept
\parskip=1.5truept
\baselineskip=11truept
\ninerm 

\item{[At]} {\bsc M.F. Atiyah}, {\ninesl Convexity and commuting
Hamiltonians}, Bull. London Math. Soc. 14 (1982), 1-15.
\item{[BM1]} {\bsc A. Banyaga and P. Molino}, {\ninesl G\'eom\'etrie des formes
de contact compl\`etement integrable de type torique}, S\'eminaire Gaston
Darboux, Montpelier, 1991-1992, 1-25.
\item{[BM2]} {\bsc A. Banyaga and P. Molino}, {\ninesl Complete Integrability
in Contact Geometry}, Pennsylvania State University preprint PM 197, 1996.
\item{[BG1]} {\bsc C. P. Boyer and  K. Galicki}, {\ninesl
3-Sasakian Manifolds}, 
to appear in {\it Essays on Einstein Manifolds}, International Press 1999;
C. LeBrun and M. Wang, Eds.
\item{[BG2]} {\bsc C. P. Boyer and  K. Galicki}, {\ninesl On Sasakian-Einstein
Geometry}, preprint DG/981108, October 1998.
\item{[BGM]} {\bsc C. P. Boyer, K. Galicki, and B. M. Mann}, {\ninesl 
The geometry and topology of 3-Sasakian manifolds},  
J. reine angew. Math. 455 (1994), 183-220. 
\item{[BGMR]} {\bsc C. P. Boyer, K. Galicki, B. M. Mann, and E. Rees}, 
{\ninesl Compact 3-Sasakian 7-manifolds with arbitrary second Betti number}, 
Invent. Math. 131 (1998), 321-344. 
\item{[Del]} {\bsc Th. Delzant}, {\ninesl Hamiltoniens P\'eriodiques et Images
Convexes de L'application Moment}, Bull. Soc. Math. France 116 (1988), 315-339.
\item{[Eli]} {\bsc Y. Eliashberg}, {\ninesl Invariants in  Contact Topology},
Doc. Math. J. DMV Extra Volume ICM II (1998), 327-338.
\item{[deMT]} {\bsc S. Falcao de Moraes and C. Tomei}, {\ninesl Moment Maps
on Symplectic Cones}, Pacific J. Math. 181(2) (1997), 357-375.
\item{[Gei]} {\bsc H. Geiges}, {\ninesl Construction of contact manifolds},
Math. Proc. Philos. Soc. 121 (1997), 455-464. 
\item{[GS]} {\bsc V. Guillemin and S. Sternberg}, {\ninesl Convexity
Properties of the Moment Map}, Invent. Math. 67 (1982), 491-513.
\item{[Gui]} {\bsc V. Guillemin}, {\ninesl Moment Maps and Combinatorial
Invariants of Hamiltonian $T^n$-spaces},$\phantom{xxx}$ Birkhauser, Boston,
1994. 
\item{[LM]} {\bsc P. Libermann and C.-M. Marle}, {\ninesl Symplectic
Geometry and Analytical Mechanics}, D. Reidel Publishing Co., Dordrecht, 1987.
\item{[LT]} {\bsc E. Lerman and S. Tolman}, {\ninesl Hamiltonian Torus Actions
on Symplectic Orbifolds and Toric Varieties}, Trans.A.M.S. 349(10) (1997),
4201-4230.
\item{[MS]} {\bsc D. McDuff and D. Salamon}, {\ninesl Introduction to
Symplectic Topology}, Oxford University Press, New York, 1995.
\item{[Tho]} {\bsc C.B. Thomas}, {\ninesl Almost regular contact manifolds}, 
J. Differential Geom. 11 (1976), 521-533. 
\item{[YK]} {\bsc K. Yano and M. Kon}, {\ninesl 
Structures on Manifolds}, Series in Pure Mathematics 3,  
World Scientific Pub. Co., Singapore, 1984. 
\item{[Zie]} {\bsc G.M. Ziegler}, {\ninesl Lectures on Polytopes},
Springer-Verlag, New York, 1995.
\medskip
\bigskip \line{ Department of Mathematics and Statistics
\hfil July 1999} \line{ University of New Mexico \hfil Revised September 1999}
\line{ Albuquerque, NM 87131 \hfil } \line{ email: cboyer@math.unm.edu,
galicki@math.unm.edu\hfil} \line{ web pages:
http://www.math.unm.edu/$\tilde{\phantom{o}}$cboyer, 
http://www.math.unm.edu/$\tilde{\phantom{o}}$galicki \hfil}

\end